\documentclass[12pt]{amsart}
\usepackage{amssymb}
\usepackage[dvips]{graphics}
\ifx\mathrm\undefined\let\mathrm\rm\fi
\ifx\mathbf\undefined\let\mathbf\bf\fi
\ifx\mathfrak\undefined\let\mathfrak\frak\fi
\ifx\mathcal\undefined\let\mathcal\cal\fi
\ifx\mathbb\undefined\let\mathbb\Bbb\fi
\ifx\emph\undefined\let\emph\it\fi
 at9.98pt

\newcommand{\n}{{{\mathfrak n}}}
\newcommand{\h}{{{\mathfrak h\,}}}

\newcommand{\Z}{{\mathbb Z}}

\newcommand{\C}{{\mathbb C}}

\newcommand{\Ref}[1]{{(\ref{#1})}}

\newcommand{\la}{\lambda}

\newcommand{\dontprint}[1]

\newcommand{\nc}{\newcommand}

\newcommand{\bs}{\boldsymbol}
\nc{\Wr}{{ {\rm Wr}}}
\newcommand{\beq}{\begin{equation}}
\newcommand{\eeq}{\end{equation}}
\newcommand{\sing}{{\rm Sing}\,}

\newcommand{\bean}{\begin{eqnarray}}
\newcommand{\eean}{\end{eqnarray}}
\newcommand{\be}{\begin{displaymath}}
\newcommand{\ee}{\end{displaymath}}
\newcommand{\bea}{\begin{eqnarray*}}
\newcommand{\eea}{\end{eqnarray*}}

\def\slg{\mathfrak{sl}}
\def\glg{\mathfrak{gl}}

{\relax}

\newtheorem%
{thm}{Theorem}[section]
\newtheorem%
{proposition}[thm]{Proposition}
\newtheorem%
{lemma}[thm]{Lemma}
\newtheorem%
{lemmadef}[thm]{Lemma-Definition}
\newtheorem%
{corollary}[thm]{Corollary}
\newtheorem%
{conjecture}[thm]{Conjecture}

\nc{\al}{\alpha}
\nc{\om}{\omega}
\nc{\La}{\Lambda}
\nc{\un}{U(\n_-)}
\nc{\px}{\frac d {d x}}

\begin{document}
\title[Higher Lame Equations and Critical Points]
{ Higher Lame Equations and Critical Points\\ of Master Functions}

\author[E. Mukhin, V. Tarasov, and A. Varchenko]
{E. Mukhin ${}^{*}$, V. Tarasov ${}^{*,\star,1}$,
\and A. Varchenko {${}^{**,2}$} }
\thanks{${}^1$\ Supported in part by RFFI grant 05-01-00922}
\thanks{${}^2$\ Supported in part by NSF grant DMS-0244579}

\maketitle

\centerline{\it ${}^*$Department of Mathematical Sciences,
Indiana University -- Purdue University,}
\centerline{\it Indianapolis, 402 North Blackford St, Indianapolis,
IN 46202-3216, USA}
\smallskip
\centerline{\it $^\star$St.\,Petersburg Branch of Steklov Mathematical
Institute}
\centerline{\it Fontanka 27, St.\,Petersburg, 191023, Russia}
\smallskip
\centerline{\it ${}^{**}$Department of Mathematics, University of
North Carolina at Chapel Hill,} \centerline{\it Chapel Hill, NC
27599-3250, USA} \medskip


\thispagestyle{empty}

\begin{abstract}
Under certain conditions, we give an estimate from above on the number
of differential equations of order $r+1$ with prescribed regular
singular points, prescribed exponents at singular points, and having a
quasi-polynomial flag of solutions. The estimate is given in terms of
a suitable weight subspace of the tensor power $U(\n_-)^{\otimes
(n-1)}$, where $n$ is the number of singular points in $\C$ and
$U(\n_-)$ is the enveloping algebra of the nilpotent subalgebra of
$\glg_{r+1}$.

\end{abstract}

\bigskip

\begin{center}
{\it Dedicated to Askold Khovanskii on the occasion
of his $60^{\text{th}}$ birthday.}
\end{center}

\section{Introduction}

Consider the differential equation
\bean\label{Ee}
F(x)\, u''(x)\ +\ G(x)\, u'(x)\ +\ H(x)\, u(x)\ =\ 0\,,
\eean
where $F(x)$ is a polynomial of degree $n$, and $G(x)\,,\ H(x)$
are polynomials of degree not greater than $n-1\,,\ n-2\,,$ respectively.
If $F(x)$ has no multiple roots, then all singular points of the equation
are regular singular. Write
\begin{eqnarray}\label{M}
F(x)\ = \ \prod_{s=1}^n \,(x-z_s)\,,
\qquad
\frac{G(x)}{F(x)}\ =\ -\ \sum_{s=1}^n\, \frac{m_s}{x-z_s}\,
\end{eqnarray}
for suitable complex numbers $m_s\,, z_s\,$. Then $0\,$ and $m_s+1\,$
are exponents at $z_s$ of equation (\ref{Ee}).
If $- l$ is one of the exponents at $\infty\,$, then
the other is $l - 1 - \sum_{s=1}^n m_s$.

\medskip\noindent
{\bf Problem} (\cite{Sz}, Ch.~6.8)\ \
{\it Given polynomials $F(x)\,,\ G(x)$ as above and a non-negative integer $l$,
\begin{itemize}
\item[(a)] find a polynomial $H(x)$
of degree at most $n-2$ such that equation (\ref{Ee})
has a polynomial solution of degree $l\,$;
\item[(b)]
find the number of solutions to Problem (a).
\end{itemize}}

If $H(x), u(x)$ is a solution to Problem (a), then the corresponding
equation \Ref{Ee} is called {\it a Lame equation} and the polynomial
$u(x)$ is called {\it a Lame function}.

\bigskip

{\bf Example.} Let $F(x) = 1 - x^2$, $G(x) = \alpha - \beta + (\alpha
+ \beta + 2) x$. Then $H = l ( l + \alpha + \beta + 1)$ and the
corresponding polynomial solution of degree $l$, normalized by the
condition $u(1) = \binom{l+\alpha}l$, is called {\it the Jacobi
polynomial} and denoted by $P^{(\alpha, \beta)}_l(x)$. \index{Jacobi
polynomial}

\bigskip

The following result is classical.

\medskip

\begin{thm}[\rm Cf.~\cite{Sz}, Ch.~6.8, \cite{St}]
\label{PS}
${}$

\begin{itemize}
\item
Let $u(x)$ be a polynomial
solution of (\ref{Ee}) of degree $l$ with roots $t_1^0,\dots, t_l^0\,$
of multiplicity one. Then $\bs t^0=(t_1^0,\dots,t_l^0)\,$ is a critical
point of the function
\bea
\Phi_{l,n}(\bs t;\bs z; \bs m)\
=\ \prod_{i=1}^{l}\prod_{s=1}^n
(t_i-z_s)^{-m_s} \prod_{1\leq i<j\leq l}
(t_i-t_j)^{2}\ ,
\eea
where $\bs z = (z_1,\dots,z_n)$ and $\bs m = (m_1, \dots , m_n)\,$.
\item
Let $\bs t^0$ be a critical point of the function
$\Phi_{l,n}(\,\cdot\,; \bs z; \bs m)\,$, then the polynomial $u(x)$ of
degree $l$ with roots $t_1^0,\dots, t_l^0$ is a solution of (\ref{Ee})
with $H(x)=(-F(x)u''(x)-G(x)u'(x))/u(x)$ being a polynomial of degree
at most $n-2\,$.
\end{itemize}
\end{thm}

The function $\Phi_{l,n}(\,\cdot\,; \bs z; \bs m)\,$ is called {\it
the master function}. The master function is a symmetric function of
the variables $t_1,\dots, t_l$. Therefore, the symmetric group $S_l$
naturally acts on the set of critical points of the master function by
permuting the coordinates.

\medskip

By Theorem \ref{PS}, the $S_l$-orbits of critical points are in
one-to-one correspondence with solution $H(x), u(x)$ of Problem (a)
such that $u(x)$ has no multiple roots.

\medskip

The following result is also classical.

\begin{thm}
[ Cf. \cite{Sz}, Ch.~6.8, \cite{H}, \cite{St} ]
\label{thm classical}
If $z_1,\dots,z_n$ are distinct real numbers and $m_1,\dots,m_n$ are negative
numbers,
then the number of solutions to Problem (a) is equal to $\binom{l+n-2}l$.
\end{thm}

Under these conditions on $\bs z$ and $\bs m$, the master function has
exactly $\binom{l+n-2}l$ $S_l$-orbits of critical points, see \cite
{Sz} and \cite{V3}.

\bigskip

The number $\binom{l+n-2}l$ has the following representation
theoretical interpretation. The universal enveloping algebra
$U(\n_-)$ of the nilpotent subalgebra $\n_- \subset \glg_2$ is
generated by one element $e_{21}$ and is weighted by powers of the
generator. The number $\binom{l+n-2}l$ is the dimension of the
weight $l$ part of the tensor power $U(\n_-)^{\otimes (n-1)}$.

\bigskip

The case of nonnegative integers $m_1,\dots,m_n$ is interesting for
applications to the Bethe ansatz method in the Gaudin model. For a
nonnegative integer $m$, let $L_m$ be the $m+1$-dimensional
irreducible $\slg_2$-module. If $m_1,\dots,m_n$ are nonnegative
integers and $m_\infty = \sum_{s=1}^n m_s -2l$ is a nonnegative
integer, then for any distinct $z_1,\dots,z_n$ the number of
$S_l$-orbits of critical points of the master function function
$\Phi_{l,n}(\,\cdot\,; \bs z; \bs m)\,$ is not greater than the
multiplicity of the $\slg_2$-module $L_{m_\infty}$ in the tensor
product $\otimes_{s=1}^n L_{m_s}$. Moreover, for generic
$z_1,\dots,z_n$, the number of $S_l$-orbits is equal
to that multiplicity, see \cite{ScV}.

\bigskip

The goal of this paper is to generalize the formulated results. We
consider linear differential equations of order $r+1$ with regular
singular points only, located at $z_1,\dots,z_n,\infty$, and
having prescribed exponents at each of the singular points. We
introduce the notion of a quasi-polynomial flag of solutions for such
a differential equation. (For a second order differential equation, a
quasi-polynomial flag is just a solution of the form
$y(x)\prod_{s=1}^n (x-z_s)^{\la_s}$ where $y(x)$ is a polynomial and
$\la_1,\dots,\la_n$ are some complex numbers.)\
We prove two facts:
\begin{itemize}
\item
Differential equations with a quasi-polynomial flag, prescribed
singular points and exponents are in one-to-one correspondence with
suitable orbits of critical points of some $\glg_{r+1}$-master
function.

\item
Under explicit generic conditions on exponents (they must be {\it
separated}), we show that the number of orbits of critical points of
the corresponding master function is not greater than the dimension of
a suitable weight subspace of $U(\n_-)^{\otimes (n-1)}$, where $\n_-$
is the nilpotent subalgebra of $\glg_{r+1}$.
\end{itemize}
The proofs are based on results from \cite{MV2, BMV}.

\bigskip

This paper was motivated by discussions with B. Shapiro of his
preprint \cite{BBS} in which another generalization of Problems (a)
and (b) is introduced for differential equations of order $r+1$. The
authors thank B. Shapiro for useful discussions of his preprint.

\section{Critical points of $\glg_{r+1}$-master functions}
\label{sec master}

\subsection{Lie algebra $\glg_{r+1}$}
\label{sec Lie algebra glg}

Consider the Lie algebra $\glg_{r+1}$ with
standard generators $e_{a,b}$, \ $a,b=1,\dots,r+1$, and
Cartan decomposition
$\glg_{r+1}=\n_-\oplus \h\oplus \n_+$,
\bea
\n_- = \oplus_{a>b} \C\cdot e_{a,b}\ ,
\qquad
\h = \oplus_{a=1}^{r+1} \C\cdot e_{a,a}\ ,
\qquad
\n_+ = \oplus_{a<b} \C\cdot e_{a,b}\ .
\eea
Set $h_a=e_{a,a} - e_{a+1,a+1}$ for $a=1,\dots,r$.
Let
$e^*_{1,1}, \dots, e^*_{r+1,r+1} \in \h^*$ be the basis dual to
the basis $e_{1,1}, \dots, e_{r+1,r+1} \in \h$.
Set $\al_a=e^*_{a,a} - e^*_{a+1,a+1}$ for $a=1,\dots,r$.
Fix the scalar product on $\h^*$ such that $(e^*_{a,a},e^*_{b,b}) = \delta_{a,b}$.

Let $\un$ be the universal enveloping algebra of $\n_-$,
\bea
\un\ =\ \oplus_{\bs l \in \Z_{\geq 0}^r}\ \un [\bs l]\ ,
\eea
where for $\bs l = (l_1,\dots,l_r)$,\ the space $\un [\bs l]$
consists of elements $f$ such that
$$
[f,h]\ = \ \langle h\,,\, \sum_{i=1}^r\,l_i\al_i \rangle f\ .
$$
The element $\prod_i e_{a_i,b_i}$ with $a_i>b_i$
belongs to the graded subspace $\un\,[\bs l]$, where
$\bs l = \sum_i \bs l_i$ with 
$$
\bs l_i\ =\ (0,0, \dots, 0, 1_b,1_{b+1},
\dots, 1_{a-1}, 0,0, \dots, 0)\ .
$$

Choose an order on the set of
elements $e_{a,b}$ with $r+1 \geq a>b\geq 1$. Then the ordered products
$\prod_{a>b} e_{a,b}^{n_{a,b}}$ form a graded basis of $\un$.

The grading of $\un$ induces the grading of $\un^{\otimes k}$ for any
positive integer $k$, \ 
$$
\un^{\otimes k}\, =\, \oplus_{\bs l}\, \un^{\otimes k}[\bs l]\ .
$$
Denote
\bea
d (k, \bs l) \ =\ \dim\ \un^{\otimes k}[\bs l]\ .
\eea

\medskip

For a weight $\La \in \h^*$, denote by $L_\La$ the irreducible
$\glg_{r+1}$-module with highest weight $\La$. Let 
$$
\bs \La =
(\La_1,\dots,\La_n) \ ,
\qquad 
\La_s \in \h^*\ ,
$$ 
be a collection of weights and
$\bs l = (l_1,\dots,l_r)$ a collection of nonnegative integers. Let
$$
L_{\bs \La}\ =\ L_{\La_1} \otimes \dots \otimes L_{\La_n}
$$ 
be the tensor
product of irreducible $\glg_{r+1}$-modules. Let $L_{\bs \La}
\,=\,\oplus_{\bs l\in \Z_{\geq 0}^r}\,L_{\bs \La}[\bs l]$ be its weight
decomposition, where $L_{\bs \La}[(l_1, \dots,l_r)]$ is the subspace
of vectors of weight $\sum_{s=1}^n\La_s - \sum_{i=1}^r l_i\al_i$. Let
$\sing L_{\bs \La}[ \bs l] \subset L_{\bs \La}[\bs l]$ be the subspace
of singular vectors, i.e. the subspace of vectors annihilated by $\n_+$. Denote
\bea
\delta (\bs \La, \bs l) \ =\ \dim\ \sing L_{\bs \La}[ \bs l] \ .
\eea
It is well-known that for given $\bs l$ and a generic set of weights
$\bs \La$, we have
\bea
d (n-1, \bs l) \ =\ \delta (\bs \La, \bs l)\ .
\eea

For given $\bs \La$ and $\bs l$, introduce the weight
\bea
\La_\infty\ =\ \sum_{s=1}^n \, \La_s\,-\,\sum_{i=1}^r\, l_i \al_i \ \in \ \h^*
\eea
and the sequences of numbers
$\bs m_1, \dots , \bs m_{r+1},
\bs m_\infty$, where
\bea
\bs m_{s} \ =\ \{\, m_{s,1} , \dots , m_{s,r+1}\,\}
\qquad
{\rm and}
\qquad
m_{s,i} \ =\ \langle \,\La_s\,,\, e_{i,i}\,\rangle\ .
\eea
Having sequences $\bs m_1, \dots , \bs m_{r+1},
\bs m_\infty$, we can recover $\bs \La$ and $\bs l$ as follows:
\bean
\label{admissible}
\La_s \ =\ \sum_{i=1}^{r+1}\ m_{s,i}\, e^*_{i,i}
\qquad
{\rm and}
\qquad
\sum_{i=1}^r \ l_i\,\al_i\ =\ \sum_{s=1}^n\ \La_s\ -\ \La_\infty
\ .
\eean
We say that the set of
sequences $\bs m_1, \dots , \bs m_{r+1}, \bs m_\infty$
of complex numbers
is {\it admissible}, if all of
the numbers $\bs l = (l_1,\dots,l_r)$,
defined by \Ref{admissible}, are nonnegative integers.

\subsection{Master functions}
\label{sec master functions}
Let $\bs z = (z_1,\dots, z_n) \in \C^n$ be a point with distinct coordinates.
Let 
$$
\bs \La\ =\ (\La_1,\dots,\La_n)\ , 
\qquad
\La_s \in \h^* \ ,
$$ 
be a collection of
$\glg_{r+1}$-weights
and $\bs l = (l_1,\dots,l_r)$ a collection of nonnegative integers.
Set $l = l_1+\dots+l_r$.
Introduce a function of $l$ variables
\bea
&
\bs t = (t^{(1)}_{1},\dots,t_{l_1}^{(1)},\dots,
t^{(r)}_{1},\dots,t_{l_{r}}^{(r)})\
\eea
by the formula
\bean
\label{master}
&&
\Phi(\bs t;\bs z;\bs \La; \bs l) =
\\
&&
\phantom{aaa}
=\ \prod_{i=1}^r \prod_{j=1}^{l_i}\prod_{s=1}^n
(t_j^{(i)}-z_s)^{-(\La_s, \al_i)} \prod_{i=1}^r\prod_{1\leq j<s\leq
l_i} (t_j^{(i)}-t_s^{(i)})^{2}
\prod_{i=1}^{r-1}\prod_{j=1}^{l_i}\prod_{k=1}^{l_{i+1}}
(t_j^{(i)}-t_k^{(j+1)})^{-1} \ .
\notag
\eean
The function $\Phi$ is a (multi-valued) function of $\bs t$, depending on
parameters $\bs z, \bs \La$. The function is called {\it the master function}.

The master functions arise in the hypergeometric solutions
of the KZ equations
\cite{SV, V1} and in the Bethe ansatz method for the Gaudin
model \cite{RV, ScV, MV1, MV2, MV3, V2, MTV}.

The product of symmetric groups 
$$
S_{\bs l}\ =\ S_{l_1}\times \dots \times
S_{l_r}
$$ 
acts on variables $\bs t$ by permuting the coordinates with the
same upper index. The master function is $S_{\bs l}$-invariant.

A point $\bs t$ with complex coordinates is called {\it a critical
point} of $\Phi(\,\cdot\,;\bs z; \bs \La; \bs l)$ if the following system of
$l$ equations is satisfied
\bean\label{BAE}
\sum_{s=1}^{n}
\frac{(\La_s,\al_1)}{t_j^{(1)}-z_s}
-\sum_{s=1,\ s\neq j}^{l_1}
\frac{2}{t_j^{(1)}-t_s^{(1)}}+\sum_{s=1}^{l_2}
\frac{1}{t_j^{(1)}-t_s^{(2)}}=0\ ,
\\
\sum_{s=1}^{n}
\frac{(\La_s,\al_i)}{t_j^{(i)}-z_s}
-\sum_{s=1,\ s\neq j}^{l_i}\frac{2}{t_j^{(i)}-t_{s}^{(i)}}+\sum_{s=1}^{l_{i-1}}
\frac1{t_j^{(i)}-t_s^{(i-1)}}+\sum_{s=1}^{l_{i+1}}
\frac{1}{t_j^{(i)}-t_s^{(i+1)}}=0\ ,
\notag
\\
\sum_{s=1}^{n}
\frac{(\La_s,\al_r)}{t_j^{(r)}-z_s}
-\sum_{s=1,\ s\neq j}^{l_{r}}
\frac{2}{t_j^{(r)}-t_s^{(r)}}+\sum_{s=1}^{l_{r-1}}
\frac{1}{t_j^{(r)}-t_s^{(r-1)}}=0\ ,
\notag
\eean
where $j=1,\dots,l_1$ in the first group of equations,
$i=2,\dots,r-1$ and $j=1,\dots,l_i$ in the second group of
equations, $j=1,\dots,l_{r}$ in the last group of equations.

In other words, a point $\bs t$ is a critical point if
\bea
\left(\Phi^{-1}
\frac{\partial \Phi }{\partial t_j^{(i)}}\right)(\bs t; \bs z; \bs \La; \bs l)\ =\ 0\ ,
\qquad
i = 1 , \dots , r,\ j = 1 , \dots l_i\ .
\eea
In the Gaudin model, equations \Ref{BAE} are called {\it the Bethe
ansatz equations}.

The set of critical points is $S_{\bs l}$-invariant.

\subsection{The case of isolated critical points}
\label{sec isolated}
We say that the pair $\bs \La, \bs l$ is {\it separating} if
\bea
(2\La_\infty\ + \ \sum_{i=1}^r c_i\al_i,\ \sum_{i=1}^r c_i\al_i)\
+\
2 \sum_{i=1}^r c_i \
\neq\ 0\
\eea
for all sets of integers $\{c_1, \dots, c_r\}$ such that\
$0\leq c_i\leq l_i$,\ {} $\sum_i c_i\neq 0$.

For example, if $\La_\infty$ is dominant integral, then
$\bs \La, \bs l$ is separating.

\begin{lemma}[Theorem 16 \cite{ScV}, Lemma 2.1 \cite{MV2}]
\label{lem finiteness}
If the pair ${\bs\La}, \bs l$ is separating, then the set of critical
points of the master function $\Phi(\,\cdot\,;\bs z;\bs \La; \bs l)$
is finite.
\end{lemma}

By Lemma \ref{lem finiteness},
for given $\bs l$ and generic $\bs \La$, the master function
$\Phi(\,\cdot\,;\bs z;\bs \La; \bs l)$ has finitely many critical points.

\begin{thm}
\label{thm crit point estimate}
Assume that the pair ${\bs\La}, \bs l$ is separating, then the number
of $S_{\bs l}$-orbits of critical points
counted with multiplicities is not greater than $d(n-1, \bs l)$.
\end{thm}

The multiplicity of a critical point $\bs t$ is the multiplicity of $\bs t$ as a solution of
system \Ref{BAE}.

\begin{proof} It is shown in \cite{BMV} that, if
$\La_1,\dots, \La_n, \La_\infty$ is a collection of dominant integral weights, then
the number of $S_{\bs l}$-orbits of critical points
counted with multiplicities is not greater than $\delta(\bs \La, \bs l)$. Together with
equality $d (n-1, \bs l) = d (\bs \La, \bs l)$ for generic $\bs \La$,
this proves the theorem.
\end{proof}

\section{Differential operators with
quasi-polynomial flags of solutions}

\subsection{Fundamental differential operator}
\label{Fundamental differential operator}
For the $S_{\bs l}$-orbit of a critical point $\bs t$ of the master function
$\Phi(\,\cdot\,;\bs z;\bs \La; \bs l)$,
define the tuple $\bs y^{\bs t} = (y_1,\dots , y_{r})$
of polynomials in
\linebreak
variable $x$,
\bea
y_i(x) \ =\ \prod_{j=1}^{l_i}\ (x-t^{(i)}_j)\ ,
\qquad
i = 1, \dots , r\ .
\eea
Since $\bs t$ is a critical point, all fractions in \Ref{BAE} are well-defined.
Therefore, the tuple $\bs y^{\bs t}$ has the following properties:
\bean
\label{i}
&&
{\rm Every \ polynomial\ } y_i\ {\rm has\ no\ multiple\ roots.}
\\
\label{ii}
&&
{\rm Every \ pair\ of\ polynomials}\ y_i\
{\rm and} \ y_{i+1}\ {\rm has\ no\ common\ roots.}
\\
\label{iii}
&&
{\rm For \ every} \ i=1,\dots,r \ {\rm and}\ s=1,\dots,n,\ {\rm if}\
m_{s,i} - m_{s,i+1} \neq 0 , \ {\rm then}\ y_i (z_s)\neq 0 .
\eean
A tuple of monic polynomials $(y_1,\dots,y_r)$ with properties \Ref{i}-\Ref{iii}
will be called
\linebreak
{\it off-diagonal}.

Define quasi-polynomials $T_1,\dots,T_{r+1}$ in $x$ by the formula
\bean
\label{formula for T}
T_i(x) \ =\ \prod_{s=1}^n\ (x-z_s)^{- m_{s,i}}\ .
\eean
Consider the linear differential operator of order $r+1$,
\bea
\label{formula fundamental operator}
D_{ \bs t} = ( \frac{d}{dx} -
\ln' ( \frac { T_{r+1} } { y_{r} } ) )\
( \frac{d}{dx} - \ln' ( \frac {y_{r}T_{r} } {y_{r-1} } ) )
\dots ( \frac{d}{dx} - \ln' ( \frac {y_2 T_2}{ y_1 } ) ) \
( \frac{d}{dx} - \ln' ( y_1 T_1 ) ) \
\notag
\eea
where $\ln'(f)$ denotes $({df}/{dx})/f$ for any $f$. We say that $D_{\bs t}$
is {\it the fundamental operator} of the critical point $\bs t$.

\begin{thm}
\label{thm on sing of fund operator}
${}$
Assume that the pair $\bs \La,\, \bs l$ is separating. Then

\begin{enumerate}
\item[(i)]
All singular points of $D_{\bs t}$ are regular and lie in
$ z_1, \dots , z_n,\infty$. The exponents of $D_{\bs t}$ at $z_s$
are 
$$
\,-m_{s,1}\,,\,-m_{s,2}+1\,,\,\dots\,,\,-m_{s,r+1}+r
$$ 
for $s = 1, \dots , n$, and
the exponents of $D_{\bs t}$ at $\infty$ are
$$
\,m_{\infty,1} \,,\,m_{\infty,2}-1\,,\,\dots\,,\, m_{\infty,r+1}-r\ .
$$
\item[(ii)]
The differential equation $D_{ \bs t}\,u\,=\,0$ has solutions $u_1,\dots,u_{r+1}$ such that
$$
u_1\ =\ y_1 T_1
$$ 
and for $i=2,\dots, r+1$ we have
\bea
\Wr\,(u_1,\dots, u_i) \ = \ y_i \ T_i \ T_{i-1} \ \dots\ T_1 \
\eea
where $\Wr\,(u_1,\dots, u_i)$ denotes the Wronskian of $u_1,\dots, u_i$ and
$y_{r+1}=1$.
\end{enumerate}
\end{thm}

\begin{proof} Part (ii) follows from the presentation of $D_{\bs t}$ as a product.
To prove part (i) consider
the operator $\tilde D_{\bs t} = T_1^{-1}\cdot D_{\bs t}\cdot T_1$, the conjugate of
$ D_{\bs t}$ by the operator of multiplication by $T_1$. Then
\bea
\tilde D_{ \bs t} = ( \frac{d}{dx} -
\ln' ( \frac { T_{r+1} } { y_{r} T_1 } ) )\
( \frac{d}{dx} - \ln' ( \frac {y_{r}T_{r} } {y_{r-1}T_1 } ) )
\dots ( \frac{d}{dx} - \ln' ( \frac {y_2 T_2}{ y_1T_1 } ) ) \
( \frac{d}{dx} - \ln' ( y_1 ) ) \ .
\eea
It is enough to show that
\begin{enumerate}
\item[(iii)]
All\ singular \ points\ of \ $\tilde D_{\bs t}$\
{ are\ regular\ and\ lie\ in}\
$\{ z_1, \dots , z_n,\infty \}$ .
\item[(iv)]
{ The \ exponents\ of\ } $\tilde D_{\bs t}$\ { at} \ $z_s$\
{\rm are } \ {}
\bea
\,0\,,\,m_{s,1}-m_{s,2}+1\,,\,\dots\,,\,m_{s,1}-m_{s,r+1}+r\ ,
\eea
{ for}\ $s=1,\dots,n$,\ { and\
the \ exponents\ of\ } $\tilde D_{\bs t}$ \ {at} \ $\infty$\ { are}
\bea
\,-l_1 \,,\,-m_{\infty,1}+m_{\infty,2}-1-l_1\,,\,\dots\,,\,
-m_{\infty,1}+m_{\infty,r+1}-r-l_1 \ .
\eea
\end{enumerate}
If all highest weights in the collection
$\La_1,\dots,\La_n, \La_\infty$
are integral dominant, then statements (iii-iv)
are proved in \cite{MV2}, Section 5. Hence statements (iii-iv) hold for arbitrary
separating $\La_1,\dots,\La_n, \ \bs l$.
\end{proof}

\subsection{Quasi-polynomial flags}
\label{Quasi-polynomial flags}
Let $\bs z = (z_1,\dots, z_n) \in \C^n$ be a point with distinct coordinates.
Let
$\bs m_1, \dots , \bs m_{r+1}, \bs m_\infty$
be an admissible set of sequences of complex numbers
as in Section \ref{sec Lie algebra glg}.
Let
\bea
D\ =\ \frac {d^{r+1}}{d x^{r+1}}
\ +\ A_1(x) \frac {d^{r}}{d x^{r}} \ +\ \dots \ + \
A_r(x)\frac {d}{d x} \ +\ A_{r+1}(x)
\eea
be a differential operator with rational coefficients.

We say that $D$ {\it is associated with}
$\bs z$, $\bs m_1, \dots , \bs m_{r+1}, \bs m_\infty$, if
\begin{enumerate}
\item[$\bullet$]
All singular points of $D$ are regular and lie in
$z_1, \dots , z_n,\infty$.
\item[$\bullet$]
The exponents of $D$ at $z_s$
are $\,-m_{s,1}\,,\,-m_{s,2}+1\,,\,\dots\,,\,-m_{s,r+1}+r$ for $s = 1, \dots , n$.
\item[$\bullet$]
The exponents of $D$ at $\infty$ are
$\,m_{\infty,1} \,,\,m_{\infty,2}-1\,,\,\dots\,,\, m_{\infty,r+1}-r$.
\end{enumerate}

Define quasi-polynomials $T_1,\dots,T_{r+1}$ in $x$ by formula
\Ref{formula for T}.

We say that an operator $D$, associated with
$\bs z$, $\bs m_1, \dots , \bs m_{r+1}, \bs m_\infty$, {\it has a
quasi-polynomial flag}, if there exists a tuple
$\bs y = (y_1,\dots, y_r)$ of monic polynomials in $x$, such that
\begin{enumerate}
\item[(i)]
For $i=1,\dots,r$, \ $\deg\, y_i \,=\, l_i$.
\item[(ii)] The tuple $\bs y$ is off-diagonal in the sense of \Ref{i}-\Ref{iii}.
\item[(iii)]
The differential equation $D\,u\,=\,0$ has solutions $u_1,\dots,u_{r+1}$ such that
$$
u_1\ =\ y_1 T_1
$$ 
and for $i=2,\dots, r+1$ we have
\bea
\Wr\,(u_1,\dots, u_i) \ = \ y_i \ T_i \ T_{i-1} \ \dots\ T_1 \ ,
\eea
where $y_{r+1}=1$.
\end{enumerate}

\begin{proposition}
An operator $D$, associated with $\bs z$, $\bs m_1, \dots , \bs m_{r+1}, \bs m_\infty$,
has a quasi-polynomial flag if and only there exists a tuple
$\bs y = (y_1,\dots, y_r)$ of monic polynomials in $x$ with properties (i-ii) and
\begin{enumerate}
\item[(iv)]
The differential operator $D$ can be presented in the form
\bea
D = ( \frac{d}{dx} -
\ln' ( \frac { T_{r+1} } { y_{r} } ) )\
( \frac{d}{dx} - \ln' ( \frac {y_{r}T_{r} } {y_{r-1} } ) )
\dots ( \frac{d}{dx} - \ln' ( \frac {y_2 T_2}{ y_1 } ) ) \
( \frac{d}{dx} - \ln' ( y_1 T_1 ) ) \ .
\eea
\end{enumerate}
\end{proposition}

\begin{proof}
The equivalence of (iii) and (iv) follows from the next lemma.

\begin{lemma}
Let $u_1,\dots,u_{r+1}$ be any functions such that
\bea
\Wr\,(u_1,\dots, u_i) \ = \ y_i \ T_i \ T_{i-1} \ \dots\ T_1 \ ,
\eea
for $i=1,\dots,r+1$.
Then $u_1,\dots,u_{r+1}$ form a basis of the space of
solutions of the equation $Du=0$ with
\bea
D = ( \frac{d}{dx} -
\ln' ( \frac { T_{r+1} } { y_{r} } ) )\
( \frac{d}{dx} - \ln' ( \frac {y_{r}T_{r} } {y_{r-1} } ) )
\dots ( \frac{d}{dx} - \ln' ( \frac {y_2 T_2}{ y_1 } ) ) \
( \frac{d}{dx} - \ln' ( y_1 T_1 ) ) \ .
\eea
\end{lemma}

\begin{proof}
It is enough to prove that\ $D\,u_i\,=\,0$\ for $i=1,\dots,{r+1}$.
By induction on $i$, we obtain that
\bea
( \frac{d}{dx} - \ln' ( \frac {y_{i}T_{i} } {y_{i-1} } ) )
\dots ( \frac{d}{dx} - \ln' ( \frac {y_2 T_2}{ y_1 } ) ) \
( \frac{d}{dx} - \ln' ( y_1 T_1 ) ) \ u\ =\
\frac{W(u_1,\dots,u_i,u)}{y_i\,T_i\,T_{i-1}\,\dots\,T_1}
\eea
from which the statement follows.
\end{proof}
\end{proof}
\bigskip

Examples of differential operators with quasi-polynomial flags are given by
Theorem \ref{thm on sing of fund operator}. If
$\bs t$ is a critical point
of the master function $\Phi(\,\cdot\,;\bs z;\bs \La; \bs l)$, then by
Theorem \ref{thm on sing of fund operator}, the fundamental
differential operator $D_{\bs t}$ is associated with
$\bs z,\ \bs m_1, \dots , \bs m_{r+1}, \bs m_\infty$ and has a quasi-polynomial flag.

\bigskip

Having $\bs z,\ \bs m_1, \dots , \bs m_{r+1}, \bs m_\infty$, define $\bs \La$,
$\bs l$
by \Ref{admissible} and the master function
$\Phi(\,\cdot\,;\bs z;\bs \La; \bs l)$ by \Ref{master}.

Having a tuple $\bs y = (y_1,\dots,y_r)$, \
$y_i(x) \ =\ \prod_{j=1}^{l_i}\ (x-t^{(i)}_j)$, denote by
\bean
\label{coord of cr pt}
\bs t = (t^{(1)}_{1},\dots,t_{l_1}^{(1)},\dots,
t^{(r)}_{1},\dots,t_{l_{r}}^{(r)})\
\eean
a point in $\C^l$, whose coordinates are roots of the polynomials of $\bs y$.
The  tuple $\bs y$ uniquely determines the $S_{\bs l}$-orbit of $\bs t$.

\begin{thm}
\label{thm flag -> crit pt}
Assume that an operator $D$ is associated with
$\bs z,\ \bs m_1, \dots , \bs m_{r+1}, \bs m_\infty$ and has a quasi-polynomial flag.
Let $\bs y$ be the corresponding tuple of polynomials. Then the point $\bs t$, introduced
in \Ref{coord of cr pt}, is a critical point of the master function
$\Phi(\,\cdot\,;\bs z;\bs \La; \bs l)$.
\end{thm}

\begin{proof} Let $u_1,\dots,u_{r+1}$ be the set of
solutions of the differential equation $D\,u\,=\,0$ giving the
quasi-polynomial flag. Introduce the functions $\tilde y_1, \dots,
\tilde y_r$ by the formulas $u_2 = \tilde y_1 T_1$ and for
$i=2,\dots, r$,
\bea
\tilde y_i\, T_i\, T_{i-1}\, \dots \,T_1\ =\ \Wr\,(u_1,\dots, u_{i-1}, u_{i+1})\ .
\eea
The functions $\tilde y_1, \dots, \tilde y_r$ are multi-valued functions.
\bea
\bullet \!\!\!\!\!\!\!\!\!
\qquad
{\rm Singularities \ of\ all\ branches\ of \ all \ of\ these\ functions\
lie\ in}\ z_1, \dots, z_n\ .
\eea
It follows from the Wronskian identities of \cite{MV2}, that
\bea
\bullet \!\!\!\!\!\!\!\!\!
\qquad
\Wr\,(\tilde y_1, y_1)\, =\, \frac {T_{2}} {T_1}
\,y_{2}\ ,
\quad
\Wr\,(\tilde y_r, y_r)\, =\, \frac{T_{r+1}} {T_r}
\,y_{r-1}\ ,
\quad
\Wr\,(\tilde y_i, y_i)\, =\, \frac{T_{i+1}} {T_i}
\,y_{i-1}\,y_{i+1}
\eea
for $i=2,\dots,r-1$. These two properties
of functions $\tilde y_1, \dots, \tilde y_r$ imply equations \Ref{BAE}
for roots of polynomials $y_1,\dots, y_r$, see \cite{MV2}.
This shows that the point $\bs t$ is a critical point.
\end{proof}

The correspondence between critical points and differential operators with
quasi-polynomial flags is reflexive. If $\bs t$ is a critical point of
the master function $\Phi(\,\cdot\,;\bs z;\bs \La; \bs l)$, which corresponds
by Theorem \ref{thm flag -> crit pt}
to the differential operator $D$, then $D$ is the fundamental differential operator
of the critical point $\bs t$.

\subsection{Conclusion}
Let $\bs z = (z_1,\dots, z_n) \in \C^n$ be a point with distinct
coordinates. Let $\bs m_1, \dots , \bs m_{r+1}, \bs m_\infty$ be an
admissible set of sequences of complex numbers as in Section \ref{sec
Lie algebra glg}. Define $\La_1,\dots, \La_n,\La_\infty$, $\bs l$ by
\Ref{admissible}. Assume that $\bs l$ is fixed and $\La_1,\dots,
\La_n$ are separating (that is generic). Then the number of differential operators $D$,
associated with $\bs z,\ \bs m_1, \dots , \bs m_{r+1}, \bs m_\infty$
and having a quasi-polynomial flag,
is not greater than $d(n-1, \bs l)$, see Theorems \ref{thm crit point estimate}
and \ref{thm flag -> crit pt}.

It is interesting to note that if $\La_1,\dots, \La_n,\La_\infty$, are
dominant integral, then the number of differential operators $D$,
associated with $\bs z,\ \bs m_1, \dots , \bs m_{r+1}, \bs m_\infty$
and having a quasi-polynomial flag, is not greater than the
multiplicity $\delta(\bs \La, \bs l)$, which could be a smaller
number. See Theorem \ref{thm flag -> crit pt} and description in
\cite{MV2} of the critical points of the corresponding master
functions.

\goodbreak

\end{document}